\documentclass[12pt,a4paper]{article}
\usepackage{amssymb}
\usepackage{amsfonts}
\usepackage{amsmath}

\input{tcilatex}
\begin{document}

\title{Global Error Control in the Runge-Kutta Solution of a Hamiltonian
System using the RKQ Algorithm}
\author{J.S.C. Prentice \\
%EndAName
Department of Applied Mathematics\\
University of Johannesburg\\
South\ Africa}
\maketitle

\begin{abstract}
We study the effect of global error control in the numerical solution of
Hamiltonian systems. In particular, we apply the RKQ algorithm in the
numerical solution of a Hamiltonian system. This algorithm is designed to
provide stepwise control of both local and global error. A test problem
demonstrates the error control features of RKQ. Good results are obtained,
despite the fact that explicit Runge-Kutta methods have been used in RKQ,
rather than symplectic Runge-Kutta methods. This simply emphasizes the value
of stepwise global error control, as per the RKQ algorithm.
\end{abstract}

\section{Introduction}

Global error control in numerical solutions of initial-value problems is of
paramount importance. In this paper, we consider the effect of global error
control on the numerical solution of Hamiltonian systems. We solve a test
problem using the RKQ algorithm, which is designed to achieve stepwise
global error control. Significantly, we use explicit Runge-Kutta methods in
RKQ, as opposed to implicit symplectic Runge-Kutta methods, which are often
preferred for Hamiltonian systems.

\section{Relevant Concepts, Terminology and Notation}

\subsection{Hamiltonian Systems}

If the total energy of a physical system is given in the form of the \textit{%
Hamiltonian}%
\begin{equation*}
H\left( p,q\right) ,
\end{equation*}%
where the canonical coordinates $p$ and $q$ are the momentum and position,
respectively, then the evolution of the system is given by%
\begin{equation}
\frac{dq}{dt}=\frac{\partial H}{\partial p},\text{ \ }\frac{dp}{dt}=-\frac{%
\partial H}{\partial q}.  \label{Ham system}
\end{equation}%
The Hamiltonian $H\left( p,q\right) $ is a \textit{first integral} of (\ref%
{Ham system}), meaning that it is constant for all $t$. For an appropriate
set of initial values, (\ref{Ham system}) is an \textit{initial-value problem%
} (IVP). Since Hamiltonian systems are nonlinear, their solution is usually
obtained numerically - using, for example, a \textit{Runge-Kutta} (RK)
method.

For ease of presentation we consider here a two-dimensional Hamiltonian,
rather than the more general case%
\begin{equation*}
H\left( p_{1},\ldots ,p_{n},q_{1},\ldots ,q_{n}\right) ,
\end{equation*}%
but this restriction will not affect our discussion.

\subsection{Runge-Kutta Methods}

Runge-Kutta methods are very well-known, and the reader is referred to the
extensive literature. For our purposes, an RK method applied to (\ref{Ham
system}) has the form%
\begin{equation}
\left[ 
\begin{array}{c}
\widetilde{q}_{i+1} \\ 
\widetilde{p}_{i+1}%
\end{array}%
\right] =\left[ 
\begin{array}{c}
\widetilde{q}_{i} \\ 
\widetilde{p}_{i}%
\end{array}%
\right] +h_{i+1}\mathbf{F}\left( t_{i},\left[ 
\begin{array}{c}
\widetilde{q}_{i} \\ 
\widetilde{p}_{i}%
\end{array}%
\right] ,h_{i+1}\right) .  \label{RK for Ham}
\end{equation}%
In (\ref{RK for Ham}), $i$ denotes discrete points along the $t$-axis, so
that $\widetilde{q}_{i}$ and $\widetilde{p}_{i}$ are numerical solutions at $%
t_{i};$ the stepsize $h_{i+1}\equiv t_{i+1}-t_{i};$ and $\mathbf{F}$ is a
vector function associated to the RK method under consideration. Initial
values of $p$ and $q$ are specified at $t_{0}.$

We define the local and global errors of an RK method as%
\begin{equation*}
\text{local error: \ \ \ }\varepsilon _{i+1}\equiv \left( \left[ 
\begin{array}{c}
q_{i} \\ 
p_{i}%
\end{array}%
\right] +h_{i+1}\mathbf{F}\left( t_{i},\left[ 
\begin{array}{c}
q_{i} \\ 
p_{i}%
\end{array}%
\right] ,h_{i+1}\right) \right) -\left[ 
\begin{array}{c}
q_{i+1} \\ 
p_{i+1}%
\end{array}%
\right]
\end{equation*}%
and%
\begin{equation}
\text{global error: \ \ \ }\mathbf{\Delta }_{i+1}\equiv \left[ 
\begin{array}{c}
\delta _{q,i+1} \\ 
\delta _{p,i+1}%
\end{array}%
\right] \equiv \left[ 
\begin{array}{c}
\widetilde{q}_{i+1} \\ 
\widetilde{p}_{i+1}%
\end{array}%
\right] -\left[ 
\begin{array}{c}
q_{i+1} \\ 
p_{i+1}%
\end{array}%
\right] .  \label{global error}
\end{equation}%
Note the use of the exact values of $q$ and $p$ in the definition of the
local error.

If the RK method is of order $r$ (denoted RK$r)$, we have%
\begin{equation*}
\varepsilon _{i+1}\propto h_{i+1}^{r+1}
\end{equation*}%
and 
\begin{equation*}
\mathbf{\Delta }_{i+1}\propto h^{r},
\end{equation*}%
where $h$ is representative of the stepsizes (for example, $h$ could be
taken as the average stepsize along the discretized $t$-axis).

Usually, \textit{symplectic} RK methods are used to solve Hamiltonian
systems \cite{hairer}. Such methods have the property that the Hamiltonian
arising from the numerical solution is bounded in the sense that it does not
drift away from the exact value; rather, it exhibits small oscillations in
the vicinity of the exact value. Also, symplectic RK methods reproduce
closed orbits in the $\left( q,p\right) $ phase space, as expected when
periodic solutions are present. The symplectic property of Hamiltonian
systems refers to the invariance of the differential 2-form $dp\wedge dq,$ a
feature that is respected by symplectic RK methods (hence their name).
Practically speaking, however, it is the first two properties listed here -
essentially constant numerical Hamiltonian and closed phase space
trajectories - that are of primary interest.

The disadvantage of symplectic RK methods is that they are \textit{implicit}
(a particular characteristic of the function $\mathbf{F})$, and this
requires the solution of a nonlinear system of equations at each node $%
t_{i}, $ which can be computationally expensive. \textit{Explicit} RK
methods do not require the solution of such a nonlinear system, but neither
are they symplectic.

\subsection{The RK$rv$Q$z$ Algorithm}

We will not discuss RK$rv$Q$z$ in detail here; the reader is referred to our
previous work where the algorithm has been discussed extensively \cite%
{prentice 2},\cite{prentice 3}. It is sufficient to state that RK$rv$Q$z$
uses RK$r$ and RK$v$ to control\ the local error via so-called \textit{local
extrapolation}, while simultaneously using RK$z$ to keep track of the global
error in the RK$r$ solution (we have $z\gg r,v)$. Such global error arises
from the propagation of the RK$v$ global error in the RK$r$ method, as a
consequence of local extrapolation. RK$rv$Q$z$ is designed to estimate the
various components of the global error in RK$r$ and RK$v$ at each node $%
t_{i} $ and, when the global error exceeds a user-defined tolerance, a 
\textit{quenching} procedure is carried out. This simply involves replacing
the RK$r$ and RK$v$ solutions with the much more accurate RK$z$ solution,
whenever necessary, so that the RK$r$ and RK$v$ global errors do not
accumulate beyond the desired tolerance.

\section{The Effect of Global Error Control}

From (\ref{global error})\ we have%
\begin{eqnarray*}
\widetilde{q}_{i+1} &=&q_{i+1}+\delta _{q,i+1} \\
\widetilde{p}_{i+1} &=&p_{i+1}+\delta _{p,i+1}.
\end{eqnarray*}%
This gives%
\begin{eqnarray*}
H\left( \widetilde{p}_{i+1},\widetilde{q}_{i+1}\right)  &=&H\left(
p_{i+1}+\delta _{p,i+1},q_{i+1}+\delta _{q,i+1}\right)  \\
&=&H\left( p_{i+1},q_{i+1}\right) +\frac{\partial H\left(
p_{i+1},q_{i+1}\right) }{\partial p}\delta _{p,i+1}+\frac{\partial H\left(
p_{i+1},q_{i+1}\right) }{\partial q}\delta _{q,i+1}+\ldots  \\
&\approx &H\left( p_{i+1},q_{i+1}\right) +\frac{dq}{dt}\delta _{p,i+1}-\frac{%
dp}{dt}\delta _{q,i+1}
\end{eqnarray*}%
where we have ignored higher-order terms. Now, say $\delta $ is an upper
bound on the magnitude of the global errors. Hence,%
\begin{equation*}
\left\vert H\left( \widetilde{p}_{i+1},\widetilde{q}_{i+1}\right) -H\left(
p_{i+1},q_{i+1}\right) \right\vert \leqslant \left( \left\vert \frac{dq}{dt}%
\right\vert +\left\vert \frac{dp}{dt}\right\vert \right) \delta .
\end{equation*}%
We see that the bound on the error in the numerical Hamiltonian is
proportional to the bound $\delta .$

For a point on the trajectory in phase space, we have%
\begin{equation*}
\left\Vert \left( \widetilde{q}_{i+1},\widetilde{p}_{i+1}\right) -\left(
q_{i+1},p_{i+1}\right) \right\Vert =\left\Vert \left( \delta _{q,i+1},\delta
_{p,i+1}\right) \right\Vert \leqslant \left\Vert \left( \delta ,\delta
\right) \right\Vert
\end{equation*}%
where $\left\Vert \cdots \right\Vert $ is any norm suitable for determining
distances in the phase space. Hence, the bound on the trajectory error is
proportional to $\delta .$ For the Euclidean norm and the two-dimensional
case considered here, we have%
\begin{equation*}
\left\Vert \left( \widetilde{q}_{i+1},\widetilde{p}_{i+1}\right) -\left(
q_{i+1},p_{i+1}\right) \right\Vert \leqslant \sqrt{2}\delta .
\end{equation*}

The implication of the above analysis is obvious: if we apply RK$rv$Q$z$ to
the problem, with a suitable tolerance of $\delta ,$ then we would generate
solutions $\widetilde{p}$ and $\widetilde{q}$ for which the error bounds on
the numerical Hamiltonian and the phase-space trajectories are acceptably
small. Furthermore, this can be achieved using explicit RK methods in the RK$%
rv$Q$z$ algorithm, as opposed to the more computationally intensive
symplectic RK methods.

\section{Numerical Example}

As an example, we consider the Hamiltonian%
\begin{equation*}
H\left( p,q\right) =\frac{p^{2}}{2}-\left( 1-\frac{p}{6}\right) \cos q
\end{equation*}%
which yields%
\begin{eqnarray*}
\frac{dq}{dt} &=&p+\frac{\cos q}{6} \\
\frac{dp}{dt} &=&\left( \frac{p}{6}-1\right) \sin q.
\end{eqnarray*}%
This is the same example considered in Hairer et al \cite{hairer}. We use
initial values%
\begin{eqnarray*}
q\left( 0\right) &=&\arccos \left( -0.8\right) \\
p\left( 0\right) &=&0
\end{eqnarray*}%
and we integrate over $t\in \left[ 0,4000\right] .$

We solve the system using RK34Q8 $\left( r=3,v=4,z=8\right) ,$ with a
tolerance of $\delta =10^{-6}$ on both the local and global error. For
comparison, we solve the system via local extrapolation only with RK3 and
RK4, also subject to a tolerance (on the local error) of $\delta =10^{-6}.$
Results are shown in the figures following the references. The various
explicit RK methods used here are the same as those referenced in \cite%
{prentice 2},\cite{prentice 3}. We use the notation RK34 to indicate the
local extrapolation algorithm using RK3 and RK4.

In Figure 1, the top two plots show the first few periods of the solution
numerical $\widetilde{p}\left( t\right) $ and $\widetilde{q}\left( t\right) ,
$ demonstrating their periodic character. Consistent with this periodicity
is the closed trajectory in phase space, shown in the third plot in Figure
1. The numerical Hamiltonian is shown in the fourth plot, for RK34Q8 and
RK34. Clearly, the latter exhibits a generally monotonic drift from the
exact value of $0.8$, while the former is essentially constant with slight
oscillations. For the sake of clarity we have not shown all the data points
in this plot; there are some $96000$ nodes on $[0,4000]$ in the computation,
and we show relatively few - sufficient, nonetheless, to exhibit the salient
features of the calculation (this also holds for subsequent plots in Figures
2-4). The drift in $H$ is slight - only about $3\times 10^{-6}$ over the
entire interval of integration - but it is definite. Given that $H$ should
be invariant, however, the result obtained with RK34Q8 should be preferred.

Figure 2 shows global errors in $\widetilde{p}\left( t\right) $ and $%
\widetilde{q}\left( t\right) ,$ determined from the difference of the RK34
solution and the RK8 solution. Clearly, the errors in the RK34 algorithm
grow as the integration proceeds, achieving maximal values of $\sim 0.02.$
On the other hand, the RK34Q8 algorithm, thanks to the high-order quenching
device, gives a solution that is bounded by $\delta =10^{-6},$ as desired.
This bound is indicated by the horizontal line labelled `tolerance'. The
maximal errors in this case are $9.998\times 10^{-7}$ in $\widetilde{q}%
\left( t\right) ,$ and $8.54\times 10^{-7}$ in $\widetilde{p}\left( t\right)
.$

In Figure 3 we show the trajectory error, determined using the Euclidean
norm, for both algorithms. Clearly, the trajectory error grows for RK34,
consistent with growth of global error in the solutions seen in Figure 2.
For RK34Q8, the trajectory error is bounded by $\sqrt{2}\times 10^{-6},$ as
expected. This bound is indicated by the horizontal line labelled `upper
bound'.

For completeness, we show the estimated global error in the RK8 solution in
Figure 4 (errors for both $p$ and $q$ are shown). This has been computed
with (see \cite{prentice 2},\cite{prentice 3},\cite{prentice 1})%
\begin{eqnarray*}
\mathbf{\Delta }_{i+1} &=&\mathbf{\varepsilon }_{i+1}+\left( \mathbf{I}%
_{2}+h_{i+1}\mathbf{F}_{q,p}\right) \mathbf{\Delta }_{i} \\
&\approx &\mathbf{\varepsilon }_{i+1}+\left( \mathbf{I}_{2}+h_{i+1}\mathbf{f}%
_{q,p}\right) \mathbf{\Delta }_{i}.
\end{eqnarray*}%
Here, $\mathbf{F}_{q,p}$ is the Jacobian of the function $\mathbf{F}$ in (%
\ref{RK for Ham}), $\mathbf{f}_{q,p}$ is the Jacobian of%
\begin{equation*}
\mathbf{f}\left( q,p\right) \equiv \left[ 
\begin{array}{c}
\frac{\partial H}{\partial p} \\ 
-\frac{\partial H}{\partial q}%
\end{array}%
\right] ,
\end{equation*}%
$\mathbf{I}_{2}$ is the $2\times 2$ identity matrix, and we assume 
\begin{equation*}
\mathbf{\Delta }_{0}=\left[ 
\begin{array}{c}
0 \\ 
0%
\end{array}%
\right] .
\end{equation*}%
Also, we have used Richardson extrapolation to estimate the local error $%
\mathbf{\varepsilon }_{i+1}$ in the RK8 solution \cite{prentice 2},\cite%
{Butcher}. For a reasonably small stepsize ($\sim 0.05$ in this problem) and
a high-order RK method, we expect this to be a good estimator. The largest
error in either component of the RK8 solution was $1.86\times 10^{-12}$.
This means that any global error in the RK8 solution would not `contaminate'
the RK34Q8 error control procedure, since the tolerance of $\delta =10^{-6}$
is considerably larger than $1.86\times 10^{-12}.$ We note here that it is
feasible to use a method of higher order than RK8 to estimate the global
error in the RK8 solution, even though this would probably increase the
computational effort. At the time of writing, however, we did not have
access to such a method. Nevertheless, in this regard one might consider
transforming the problem into a second-order problem, making it suitable for
Nystr\"{o}m integration, as we considered in \cite{prentice 4}. This would
enable the use of the $\left( 10,12\right) $ embedded Nystr\"{o}m pair due
to Dormand et al \cite{dormand}.

\subsection{Relative Error Control}

In this paper, we have only considered absolute error control, since the
magnitude of the solution was $\sim 1.$ If the magnitude of the solution had
been much larger than unity, it would have been better to implement relative
error control. We will not discuss this in detail; a thorough account has
been given in \cite{prentice 3}, wherein we generalized the RKQ algorithm.
It suffices to say that at each node, for relative error control in the
problem considered here, the tolerance would have the form%
\begin{equation*}
\delta =\min \left\{ \max \left\{ \delta _{A},\delta _{R}\left\vert
q_{i}\right\vert \right\} ,\max \left\{ \delta _{A},\delta _{R}\left\vert
p_{i}\right\vert \right\} \right\} ,
\end{equation*}%
where $\delta _{A}$ and $\delta _{R}$ are user-defined (the presence of $%
\delta _{A}$ caters for those situations where $q_{i}$ and/or $p_{i}$ are
very close to zero).

\section{Conclusion}

We have explored global error control in the numerical solution of
Hamiltonian systems. A theoretical analysis shows that if the error in the
solution is bounded, then errors in quantities such as the numerical
Hamiltonian and phase space trajectories are also bounded. We have
considered the use of the RKQ algorithm, with explicit Runge-Kutta methods,
as our choice of numerical integrator, since RKQ\ is specifically designed
to achieve global error control. A test problem has demonstrated the
expected results. In addition to the bounding of the error in the solution,
we also observe an essentially constant numerical Hamiltonian and a bounded
error in the numerical trajectory in phase space. This contrasts Hamiltonian
drift and unbounded trajectory deviation normally associated with the use of
explicit Runge-Kutta methods in solving Hamiltonian systems. It is our
contention that, even though RKQ utilizes explicit Runge-Kutta methods,
stepwise global error control via RKQ leads to results with a similar
quality to those that would be obtained using symplectic methods.

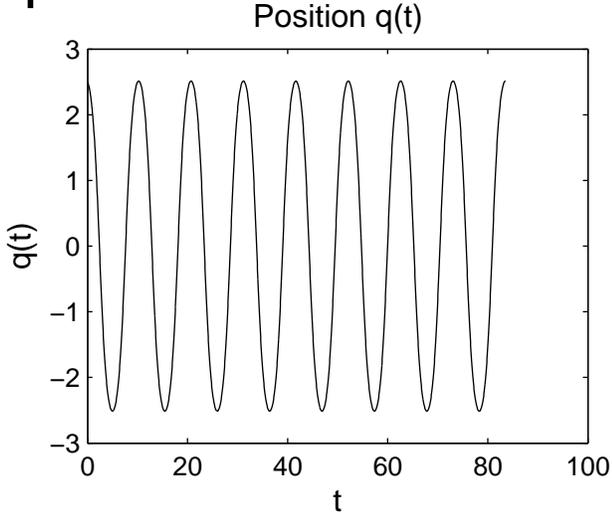
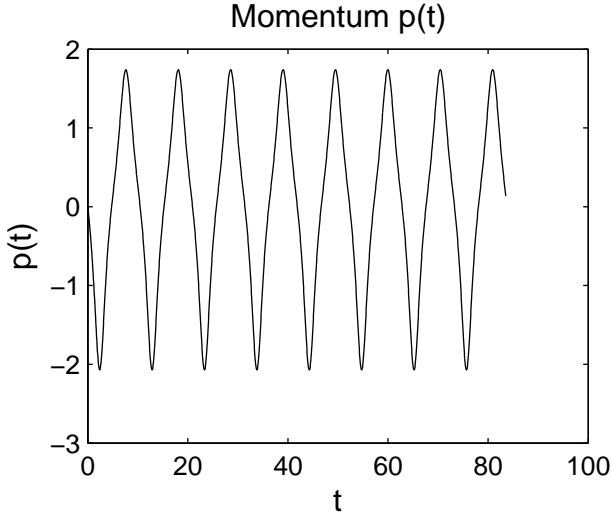
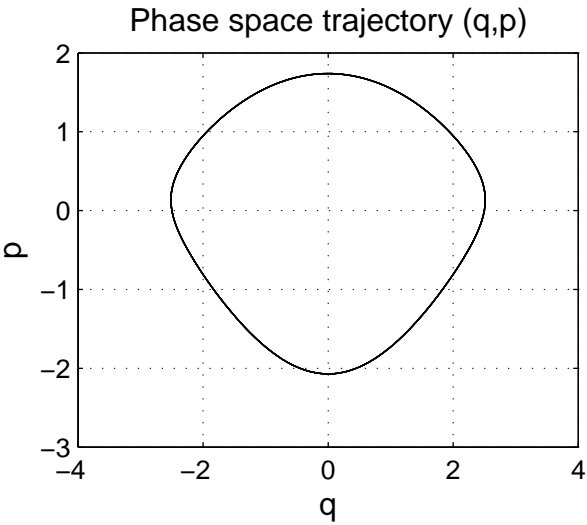
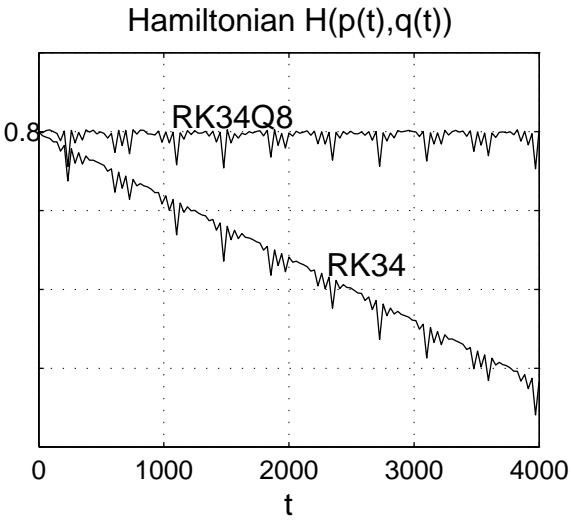
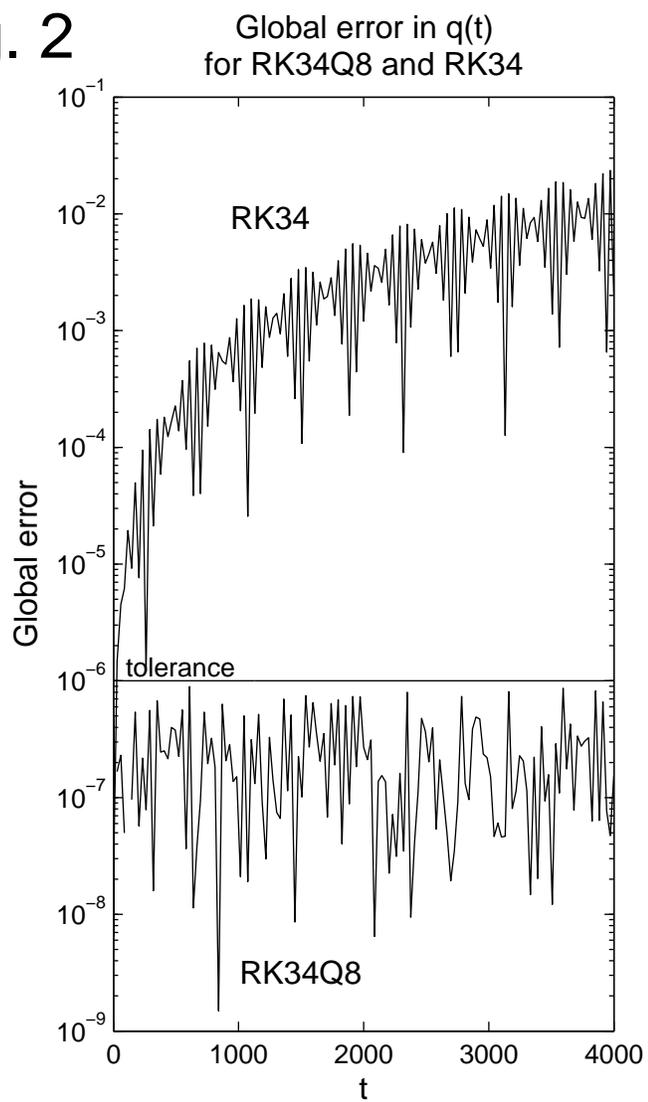
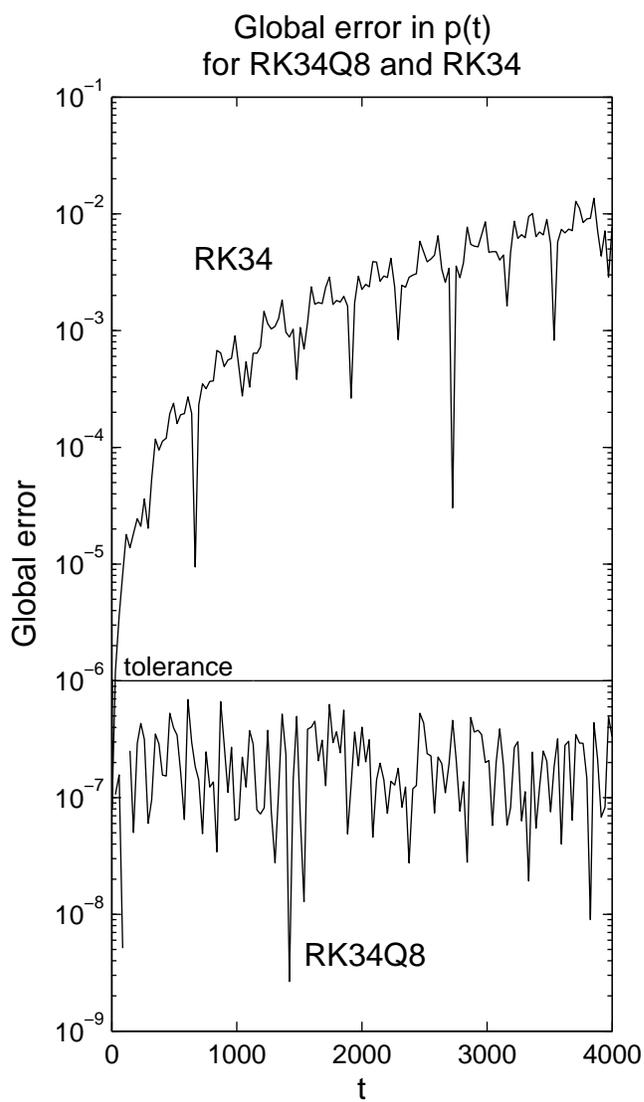
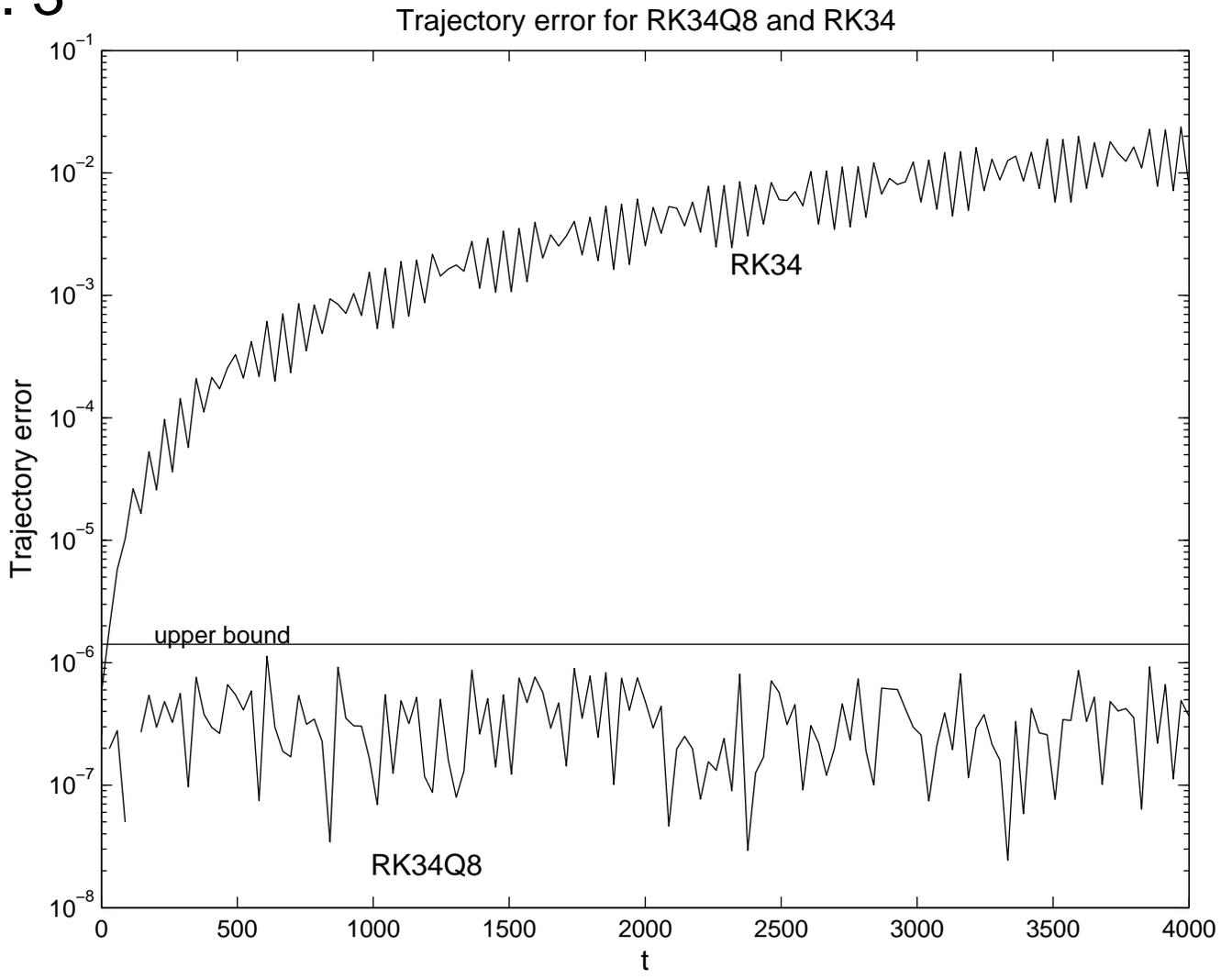
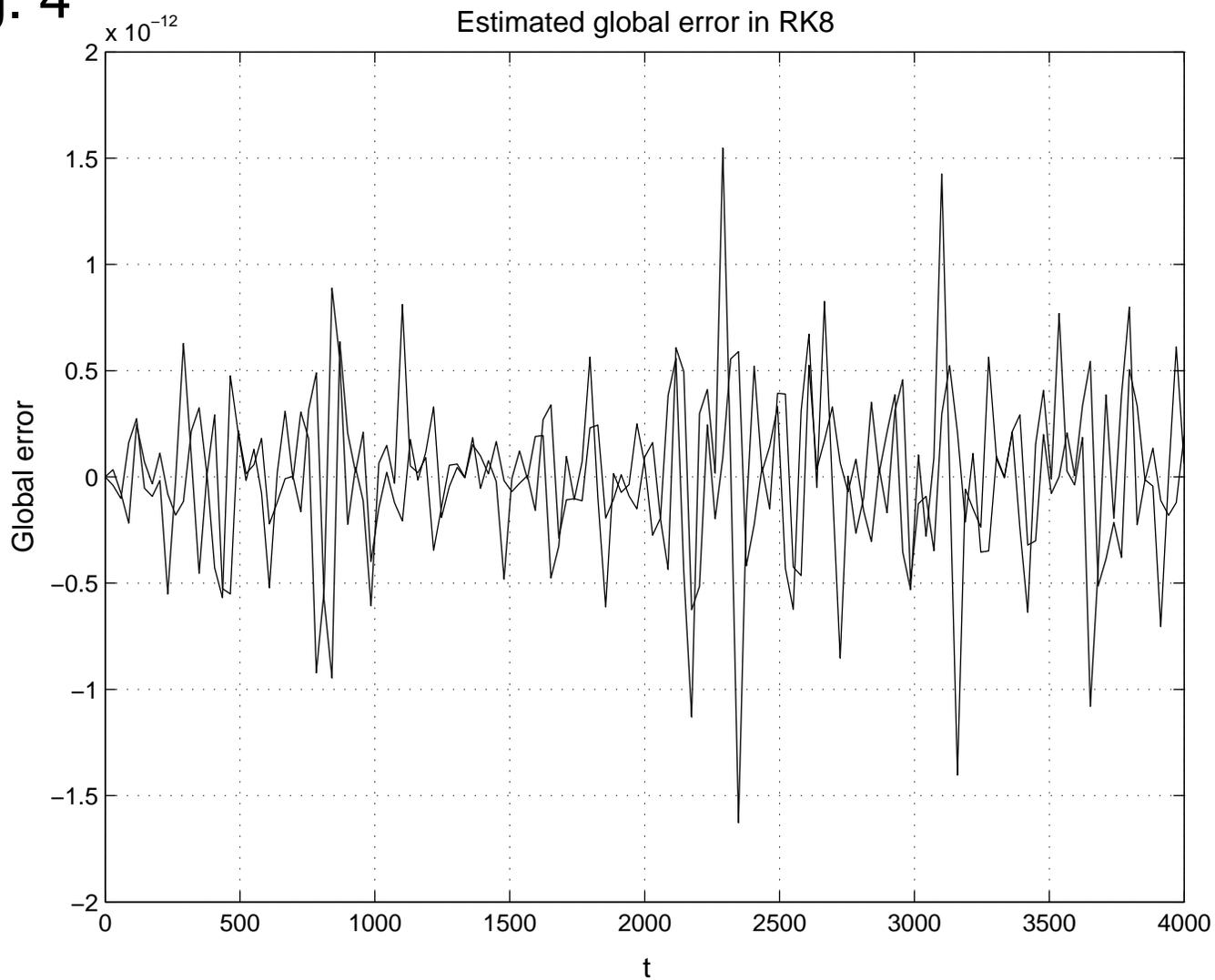
\end{document}